\documentclass[11pt,A4,twoside]{article}
\usepackage[latin1]{inputenc}
\usepackage{amssymb}
\usepackage{amsmath}
\usepackage{amsfonts}
\usepackage{amsbsy}
\usepackage{natbib}
\usepackage{eucal}
\usepackage{graphicx}
\usepackage{titlesec}
\usepackage[english]{babel}
\usepackage{times}
\usepackage{titling}
\usepackage{mathptmx}
\usepackage{float}
\usepackage[font=small]{caption}
\usepackage{array}
\usepackage[colorlinks,citecolor=blue,urlcolor=blue]{hyperref}
\usepackage{authblk}
\usepackage{amsthm}
\usepackage{bbold}

\usepackage{setspace}
\DeclareSymbolFont{txgreek}{OML}{cmr}{m}{it}

\renewcommand{\abstract}[1]{{\small\noindent
\hrulefill\par \vspace*{0.1cm}\noindent{\small\bf\sffamily
{Abstract}}\parindent=0pt\par\noindent\vspace{-0.1cm}\noindent\hrulefill\par\vspace*{0.5\baselineskip}\hspace*{0cm}\renewcommand{\baselinestretch}{1.1}\sffamily{#1}\par\vspace*{-0.1cm}\noindent\hrulefill}}

\newtheorem{theorem}{Theorem}[section]
\newtheorem{corollary}{Corollary}[theorem]

\def\and{,\;}
\DeclareMathSizes{11}{11}{7.8}{6.5}

\def\paragraf{\fontsize{9}{10pt}\fontfamily{phv}\fontshape{it}\selectfont}
\titleformat{\paragraph}
{\titlerule[0pt]
\vspace{0cm}%
\paragraf}
{\theparagraph}{.5em}{\vspace*{0\baselineskip}}

\def\titol{\fontsize{12.045}{12pt}\fontfamily{phv}\fontseries{b}\selectfont}
\titleformat{\section}
{\titlerule[0pt]
\vspace{0.2cm}%
\titol}
{\thesection.}{.5em}{\vspace*{0\baselineskip}}

\def\titolp{\fontsize{11.045}{11pt}\fontfamily{phv}\fontseries{b}\fontshape{it}\selectfont}
\titleformat{\subsection}
{\titlerule[0pt]\bigskip
\vspace{-0.4cm}%
\titolp}
{\thesubsection.}{.5em}{\vspace*{0\baselineskip}}

\def\titolpp{\fontsize{10.045}{10pt}\fontfamily{phv}\fontshape{it}\selectfont}
\titleformat{\subsubsection}
{\titlerule[0pt]
\vspace{-0.4cm}%
\titolpp}
{\thesubsubsection.}{.5em}{\vspace*{0\baselineskip}}

\usepackage{titling}
    \pretitle{\begin{center}\sffamily\fontsize{18pt}{20pt}\selectfont}
    \posttitle{\par\end{center}}
    \preauthor{\begin{center}\fontsize{12pt}{14pt}\selectfont}
    \postauthor{\par\end{center}\vspace{0bp}}
    \predate{}
    \date{}
    \postdate{}
    
\usepackage{xcolor}
\usepackage{multirow}
\def\1{{\rm l}\hskip -0.21truecm 1}
\definecolor{blau}{rgb}{0,0,0.5}

\newtheorem{definition}[theorem]{Definition}
\newtheorem{remark}[theorem]{Remark}
\newtheorem{notation}[theorem]{Notation}
\usepackage{cancel}
\DeclareMathOperator{\Tr}{Tr}

\title{Second Order Markov multistate models}

\thanksmarkseries{arabic}
\author{Mireia Besal\'u\thanks{Departament Gen\`etica, Microbiologia i Estad\'istica, Universitat de Barcelona, Barcelona, Spain} \and  Guadalupe G\'omez Melis\thanks{Departament d'Estad\'istica i Investigaci\'o operativa, Universitat Polit\`ecnica de Catalunya, Jordi Girona 1-3, 08034 Barcelona, Spain}}

\def\headers#1{\fontsize{8.5}{10}\centering\sffamily\itshape{#1}}
\def\page#1{\fontsize{8.5}{10}\sffamily{#1}}
\usepackage{fancyhdr}%Headers

\begin{document}
\maketitle

% Headers
\thispagestyle{empty}
\renewcommand{\headrulewidth}{0truecm}
\pagestyle{fancy}
\rhead[\headers{Second Order Markov multistate models}]{\page{\thepage}}
\lhead[\page{\thepage}]{\headers{M. Besal\'u, G. G\'omez Melis}}
 \lfoot{} \rfoot{}
\cfoot{}

%75-100 words

\abstract{Multistate models (MSM) %  to describe the evolution of individuals along time and through different states
are well developed for continuous and discrete times under a first order Markov assumption.
Motivated by a cohort of COVID-19 patients, an MSM was designed  based on 14 transitions among  7  states of a patient.
Since a preliminary analysis showed that the first order Markov condition was not met for some transitions, we have developed  a second order Markov model where the future evolution not only depends on the current but also on the preceding state. Under a discrete time analysis, assuming homogeneity and that past information is restricted to 2 consecutive times, we expanded the transition probability matrix and proposed an extension of the Chapman-Kolmogorov equations.
}

\paragraph{MSC: 62M09, 62N02, 60J10}

\paragraph{Keywords: Multistate models, Non-Markov, COVID-19}

\renewcommand{\baselinestretch}{1.2}
\bigskip

\section{Introduction}\label{sec1}

 Multistate models (MSM) provide a very convenient methodology  to describe  the life history of an individual which at any time occupies one of a few possible states. In particular, they are appropriate to describe the clinical course of a disease and are routinely used in research to model the progression of patients among different states.

MSM theoretical justification is  based on the theory of stochastic processes, that is, on sets of random variables representing the evolution of a process over time.  The time can be chosen to be  discrete or continuous; while discrete times assume a stepwise process where the fixed time between successive steps is not part of the model, continuous time models allow changes of the states at any time. This class of models allows for an extremely flexible approach that can model almost any kind of longitudinal failure time data. This is particularly relevant for modeling different events, which have an event-related dependence, like the occurrence of a disease changing the risk of death (\cite{Hougaard}).

The first order Markov condition establishes that the future evolution of the stochastic process only depends on the current state and is  frequently assumed in multistate models. However, this condition might often be not too realistic to describe clinical situations. To  test it \cite{TitmanPutter} develop general log-rank tests that can be applied to general multistate models under right-censoring. 

A plausible approach to lessen the first order Markov assumption is  to consider a higher order Markov process. A Markov process  of order $k$ is such that
the dependence of the process on the whole  history is only through the $k$ states previously occupied.  Although it  is often observed that higher order Markov processes can model the data better, models for Markov processes of higher order are scarcely used in practice  because they depend on a  very large number of parameters leading to computational difficulties (\cite{Ching(2003), Logan(1981)}). Most instances of higher order Markov models, which have been used so far,  involved discrete time models (known as Markov chains). \cite{Tong(1975)}  defines a $k$ order Markov chain $\{X_1, \cdots, X_n, \cdots\}$ as the one such that the conditional probabilities  satisfy 
\begin{equation}\label{MCok}
 P(X_n|X_{n-1}, X_{n-2}, \cdots) =P(X_n|X_{n-1}, X_{n-2}, \cdots, X_{n-k})  
\end{equation}
for all $n$, where $k>0$ is the smallest integer holding the above condition. 

Trying to combine realism with parsimony \cite{Raftery(1985)}  introduces a Markov chain model of order $k$ and $m$ states where \eqref{MCok} is expressed as a linear combination of contributions from each $X_{n-1}, X_{n-2}, \cdots, X_{n-k}, $ that is,
\begin{equation}\label{MCRaftery}
P(X_n=x_n|X_{n-1}=x_{n-1},  \cdots, X_{n-k}=x_{n-k}) =\sum_{i=1}^k\lambda_i q_{x_n,x_{n-i}}
\end{equation}
where $\lambda_1+\cdots +\lambda_k=1$ and $Q=\{q_{ij}\}$ is a non-negative $m\times m$ matrix with columns equal to 1 such that
$0\leq \sum_{i=1}^k\lambda_i q_{x_n,x_{n-i}}\leq 1$. 
Since the number of independent parameters is $(m-1)\times m^k$, model \eqref{MCRaftery} reduces this figure quite drastically. For instance, the  usual 100 parameters  needed for a second order process with 5 states reduces under \eqref{MCRaftery} to only 21.
 Other authors \cite{Ching(2003), Islam(2006)} have  used different relations among parameters to make the estimation  of them feasible and their  interpretation unambiguous.

In this paper we propose second order Markov multistate models as a way of enriching the pathway information and still control the number of parameters  while keeping the interpretability of the transition probabilities.
Analysis using second order Markov models are scarce. Among them, \cite{Shorrocks(1976)} investigated the Markovian assumption in modelling income mobility and concluded that transition rates should depend on both current income and immediate past history, hence,  a second order Markov model was implemented. \cite{Shamshad} uses a second order Markov model for synthetic generation of wind speed time series data.

Second order Markov models assume that the progression of the individuals not only depends on the current but also on the preceding state. Second order Markov multistate models are characterized by means of a  $M\times M \times M$ tensor, where $M$ is the number of states. In this work, we define an extended transition probability matrix as $M$ different matrices of order $M\times M$. To be able to compute $n$-step transition probabilities we extend the first order Chapman-Kolmogorov equations.
We conclude the paper  with an illustration consisting of a cohort of more than 2000 COVID-19 patients from five hospitals in the Barcelona metropolitan area who were hospitalized during the first wave of the coronavirus pandemic (March-April 2020). For this data we have built a multistate model based on 14 possible transitions among the seven states where a patient can be in after his/her admission.
Since it is seen that the first order Markov condition  does not hold for all the transitions we fit a second order Markov chain. We estimate the second order transition probabilities and based on those we compute, among others, the transition  probability from one state  to another, after a given number of hospitalized days, and differentiating between patients that arrive to the hospital with the disease from those who develop the disease at the hospital. The paper ends with a discussion on shortcomings while setting the path for future research.
 
%%%%%%%%%%%%%%%%%%%%%%%%%%%%%%%%%%%%%%%%%%%%%%
%%%%%%%%%%%%%%%%%%%%%%%%%%%%%%%%%%%%%%%%%%%%%%
\section{Characterization of first order Markov multistate processes}
%%%%%%%%%%%%%%%%%%%%%%%%%%%%%%%%%%%%%%%%%%%%%%
%%%%%%%%%%%%%%%%%%%%%%%%%%%%%%%%%%%%%%%%%%%%%%

A multistate process is a continuous-time stochastic process $X=\{X_t, \enspace t\geq 0\}$ taking values in a discrete state space ${\mathcal S}=\{1,\cdots, M\}$. We denote by ${\mathcal F}_{t}:=\sigma\{X_s:s\leq t\}$ a $\sigma$-algebra consisting on the observation of the process over the interval $[0,t]$ and we refer to it as  a filtration. We can think of a   filtration as the history of the process up to time $t$ containing  the information on the previous occupied states up until time t. 

The law of a multistate process is defined by its finite dimensional distribution and  is fully characterized through either one of the following 3 functions: transition probabilities, transition intensities or cumulative transition intensities.
The \textbf{transition probability} between states $h$ and $j$ for times $s$ and $t$, $s< t$ is defined by: 
$$P_{hj}(s,t;{\cal F}_{s-})=P(X_t=j\mid X_s=h;\,{\cal F}_{s-})\qquad \mbox{for} \qquad h,j\in {\cal  S}=\{1,\cdots, M\}$$
and denote the probability of the process  being at state $j$ at time $t$ knowing that it has  been at state $h$ at time $s$ as well as knowing all the previous trajectory before $s$.
The \textbf{transition intensity} between states $h$ and $j$, for time $t$ is defined by: 
$$\alpha_{hj}(t;\,{\cal F}_{t-})=\lim_{\Delta t\rightarrow 0} \frac{1}{\Delta t}P_{hj}(t,t+\Delta t;\, {\cal F}_{t-})\qquad \mbox{for}\qquad h,j\in {\cal S}=\{1,\cdots, M\}$$
and denotes the propensity to change from state $h$ to state $j$ at time $t$.
The \textbf{cumulative (integrated) transition intensity} between states $h$ and $j$ at time $t$ is defined by: 
$$A_{hj}(t;\,{\cal F}_{t-})=\int_0^t \alpha_{hj}(u;{\cal F}_{t-})du \qquad \mbox{for}\qquad h,j\in {\cal S}=\{1,\cdots, M\}$$

Transition probabilities, transition intensities and cumulative transition intensities are  summarized by means of  $M\times M$ matrices. In particular,  for every trajectory collected in ${\cal F}_{t-}$  and for every $s,\,t$ such that $s< t$,
we denote by $\mathbf{P}$ the transition probability matrix 
\begin{equation*}
   \mathbf{P}(s,t ;\, {\cal F}_{s-})=\{P_{hj}(s,t;\,{\cal F}_{s-}); h,j\in {\cal S}=\{1,\cdots, M\}\}.
\end{equation*}

Finally, we define the row vector of \textbf{state occupation probabilities} 
$\mathbf{\pi}(t)=(\pi_j(t))_{j\in\cal{S}}$ for each time $t$
where
 $\pi_j(t)=P(X_t=j)$ is the probability of being in state $j$ at time $t$.

%%%%%%%%%%%%%%%%%%%%%%%%%%%%%%%%%%%%%%%%%%%%%%%%%%%%%%%% 
\subsection{First order Markov and homogeneity assumptions}
%%%%%%%%%%%%%%%%%%%%%%%%%%%%%%%%%%%%%%%%%%%%%%%%%%%%%%%% 

It is clear that  some restrictions have to be made  in order to estimate the transition probabilities $P_{hj}(s,t;{\cal F}_{s-})$ for every pair of states $h$ and $j$, for every pair of times $s$ and $t$ and for all the possible  trajectories before $s$. The Markov and the homogeneity assumptions are key to make inference feasible.

\begin{definition}
A multistate process satisfies the \textbf{first order Markov assumption} if for all $ h,j\in {\cal S}=\{1,\cdots, M\}$ and $s,\,t$ such that $s<t$ 
\begin{equation*}
P_{hj}(s,t;\,{\cal F}_{s-})=P(X_t=j\mid X_s=h;\,{\cal F}_{s-}) =P(X_t=j\mid X_s=h)=P_{hj}(s,t).
\end{equation*}
That is, under the first order Markov assumption, different trajectories before $s$ will not change the transition probabilities. Under the {first order Markov assumption}, an $M\times M$ matrix, ${P}(s,t)$, is needed \textbf{for every $(s,t)$}

\begin{equation*}
\mathbf{P}(s,t)=\{P_{hj}(s,t);\; h,j\in {\cal S}=\{1,\cdots, M\}\}
\end{equation*}

\end{definition}

\begin{definition} \label{def:hom1}
A first order Markov multistate process is said to be \textbf{homogeneous} if the transition probability between any states at  given times $t,\,s \enspace (s < t)$ depends only on the difference between these two times ($t-s$), that is, 
\[P_{hj}(s,t)=P_{hj}(0,t-s)=P_{hj}(t-s).\]
In this case only a $M\times M$ matrix $\mathbf{P}(t)$ for every time $t$ is needed.
\end{definition}

%%%%%%%%%%%%%%%%%%%%%%%%%%%%%%%%%%%%%%%%%%%%%%
\subsection{Markov Test} \label{sec:MarkovTest}
%%%%%%%%%%%%%%%%%%%%%%%%%%%%%%%%%%%%%%%%%%%%%%

We should validate the Markov condition if we want to proceed analysing the data under this assumption.  \cite{TitmanPutter} discuss several tests to check the Markov assumption and propose a new one. One choice would be to include the time of entry into each state as a covariate within a   Cox model and test its significance through a likelihood ratio test (\cite{Kay}).  A second possibility would be to use the  stratified version of the Commenges-Andersen's  test to detect a shared frailty. 
Other authors (\cite{Rodriguez-Girondo}) have developed local and global tests for the  Markov conditions based upon the observed Kendall's $\tau$ for the progressive three-state illness-death model.

In this paper, and in the subsequent COVID-19 analysis, we will  validate the Markov assumption for each transition by means of \cite{TitmanPutter}'s test that we briefly describe. The main idea of this test is that under the first order Markov assumption the  rate of transitions at time $t>s$ will not be affected by the state occupied at time $s$. If we want to check the Markov assumption for the transition between the states $l$ and $m$ $(l,\,m\in\cal{S})$, we divide the subjects into two different groups: the ones that at time $s$ are in a fixed state $j\in \cal{S}$ and the ones who are not there.
Then, for each transition ($l\rightarrow m$) the null hypothesis for a fixed state $j$ and fixed  time $s$  ($s\in[t_0,t_{\max}]\subset [0,\tau]$, $\tau$ total follow-up) is stated as:
\[\mbox{H}^{(j)}_{0s}(l,m): \alpha_{lm}(t\mid X(s)=j)=\alpha_{lm}(t\mid X(s)\neq j)\quad \mbox{for any } \quad t\in[s,\tau]\]
 and can be tested with the log-rank statistic
\[U_s^{(j)}(l,m)=\sum_{i=1}^n\int_s^\tau\left\lbrace\delta_i^{(j)}(s)-\frac{\sum_{k=1}^n \delta_k^{(j)}(s) Y_{kl}(t)}{\sum_{k=1}^n Y_{kl}(t)}\right\rbrace dN_i^{(lm)}(t),\]
where $\delta_i^{(j)}(s)=\1{\{X_i(s)=j\}}$ denotes whether individual $i$ has been in state $j$ at time $s$, $Y_i(t)$ is the at risk indicator for the process $X_i(t)$, $ Y_{il}(t)=\1{\{X_i(t^-)=l\}}Y_i(t)$ is the at risk indicator of transition $l\rightarrow m$ for subject $i$ and $N_i^{(lm)}(t)$ is the counting process reporting the number of times of the transition $l\rightarrow m$ up to time $t$. 

The standardized statistics 
$$\overline{U}_s^{(j)}(l,m)=\frac{U_s^{(j)}(l,m)}{\sqrt{\widehat{Var}(U_s^{(j)})(l,m)}}.$$
can be compared to a $N(0,1).$
Moreover, \(\left\lbrace \overline{U}_s^{(j)}(l,m), s\in[t_0,t_{\max}]\right\rbrace\)
converges to a zero mean Gaussian process with a covariance function that can be consistently estimated.

Given the  null  hypothesis for a fixed state $j$
\[\mbox{H}^{(j)}_{0}(l,m): \alpha_{lm}(t\mid X(s)=j)=\alpha_{lm}(t\mid X(s)\neq j) \quad  \forall s\in[t_0,t_{\max}]\subset [0,\tau] \;  \mbox{and } \; t\in[s,\tau],\]
a global test statistic  
can be defined based on summary statistics of $\{\overline{U}_s^{(j)}(l,m), s\in[t_0,t_{\max}]\}$ such as
\(\displaystyle \int_{t_0}^{t_{\max}} \left\vert\overline{U}_s^{(j)}(l,m)\right\vert ds,\)
$\sup_{s\in[t_0,t_{\max}]} \vert\overline{U}_s^{(j)}(l,m)\vert$ or $\displaystyle\int_{t_0}^{t_{\max}} w(s)\vert\overline{U}_s^{(j)}(l,m)\vert ds$ for some weight function $w(s)$.

Finally, an overall test statistic for the null hypothesis for any possible  $j$  and   for  all  $s\in[t_0,t_{\max}]\subset [0,\tau]$  and $t\in[s,\tau]$
\[\mbox{H}_{0}(l,m)\!: \alpha_{lm}(t\mid X(s)=j)=\alpha_{lm}(t\mid X(s)\neq j) \quad \forall j,\, \forall s\in[t_0,t_{\max}]\!\subset [0,\tau]\;  \mbox{and} \; \forall t\!\in\![s,\tau]\]
can be defined from the global test statistics, for instance 
as the mean, the maximum or weighted mean of them.
These tests are implemented in \texttt{R} with the function \texttt{MarkovTest} of the package \texttt{mstate} of \citep{mstate}.
Details of how are implemented are postponed to the illustration in Section \ref{Sec: MarkovTest2}.

%%%%%%%%%%%%%%%%%%%%%%%%%%%%%%%%%%%%%%%%%%%%%%%%%%%%%%%% 
%%%%%%%%%%%%%%%%%%%%%%%%%%%%%%%%%%%%%%%%%%%%%%%%%%%%%%%% 
\section{Characterization of second order Markov multistate processes}
%%%%%%%%%%%%%%%%%%%%%%%%%%%%%%%%%%%%%%%%%%%%%%%%%%%%%%%% 
%%%%%%%%%%%%%%%%%%%%%%%%%%%%%%%%%%%%%%%%%%%%%%%%%%%%%%%% 

 Since the first order Markov assumption is an strong assumption we relax this condition defining a second order Markov assumption. The main idea is that the future evolution of the stochastic process  depends not only on the present state and but also on the preceding one.

%%%%%%%%%%%%%%%%%%%%%%%%%%%%%%%%%%%%%%%%%%%%%%%%%%%%%%%%%%%%
 \subsection{Second Order Markov Transition probabilities} \label{sec:2norderTransProb}
 %%%%%%%%%%%%%%%%%%%%%%%%%%%%%%%%%%%%%%%%%%%%%%%%%%%%%%%%%%%
 
 We start defining the second order Markov transition probabilities and describing how they can be summarized into a set of as many matrices as states.

\begin{definition}
For times $(s,t,u),\,s< t< u$ and {adjacent states} $h, j, k$ (meaning that there is a direct path between them), the probability $P_{hjk}(s,t,u;{\mathcal{F}}_{s-}) = P(X_u=k \mid X_s=h, X_t=j; \mathcal{F}_{s-})$  satisfies a second order Markov assumption if and only if
\begin{eqnarray*}
P_{hjk}(s,t,u;{\mathcal{F}}_{s-})&=&P(X_u=k \mid X_s=h, X_t=j|{\mathcal{F}}_{s-}) \\&=&P(X_u=k \mid X_s=h, X_t=j)=P_{hjk}(s,t,u).
\end{eqnarray*}
\end{definition}
\noindent
Under the second order Markov assumption the transition probabilities are summarized,  for every three times $(s,t,u),\,s<t<u$, by  an $M\times M\times M $ tensor $\mathbf{P}(s,t,u)$
  \begin{equation*}
    \mathbf{P}(s,t,u)=\{P_{hjk}(s,t,u);\, h,j,k\in {\cal S}=\{1,\cdots, M\}\}.
  \end{equation*}
In order to have a more manageable mathematical object we denote, for each state $h\in{\cal S},$ a matrix of dimension $M$, $\mathbf{P_{(h)}}(s,t,u)$ as follows:
\[\mathbf{P_{(h)}}(s,t,u)=\left(P_{hjk}(s,t,u)\right)_{j,k\in\mathcal{S}},\]
hence, %and we will able to represent
 the tensor $\mathbf{P}(s,t,u)$ of transition probabilities can be equivalently represented as $M$ matrices of order $M$ for each $s<t<u$.

\begin{remark}
The matrices $\mathbf{P_{(h)}}(s,t,u)$ are not always stochastic matrices because
\[\sum_{k\in\mathcal{S}} P_{hjk}(s,t,u)=\begin{cases}0 &\mbox{if }\; \forall j,k,  \;\mbox{the transitions} \;
h\rightarrow j \mbox{ or } j\rightarrow k  \mbox{ are not possible }\\1 &\mbox{otherwise.}\end{cases}\]
For example, if $h$ is an absorbent state all the matrices will be 0 except for the element $P_{hhh}=1$.
\end{remark}

\begin{definition} \label{def:hom2}
A second order Markov multistate process is said to be \textbf{homogeneous} if the transition probability between any three states at  given times $(s,t,u),\,s< t< u$, depends only on the differences $t-s$ and $ u-t$ between the two  consecutive times %($t-s, u-t$), 
that is, 
\[P_{hjk}(s,t,u)=P_{hjk}(t-s, u-t)\]
In this case only a $M\times M\times M$ tensor $\mathbf{P}(s,t)$ for every pair of times $(s,t)$ ($s<t$) is needed
  \begin{equation*}
    \mathbf{P}(s,t)=\{P_{hjk}(s,t);\, h,j,k\in {\cal S}=\{1,\cdots, M\}\}.
  \end{equation*}
Using the previous notation, denote as $\mathbf{P_{(h)}}(s,t)$ the matrix of dimension $M$ for each state $h\in{\cal S}$  and for every pair of times $(s,t)$ ($s<t$), that is, 
\[\mathbf{P_{(h)}}(s,t)=\left(P_{hjk}(s,t)\right)_{j,k\in\mathcal{S}}.\]
Note that the tensor $\mathbf{P}(s,t)$ of transition probabilities under homogeneity can be equivalently represented as  $M$ matrices of dimension $M$ for each two times $(s,t); s<t$ where $s$ stands for the time from $h$ to $j$ and $t$ stands for the time from $j$ to $k$.
\end{definition}

%%%%%%%%%%%%%%%%%%%%%%%%%%%%%%%%%%%%%%%%%%%%%%%%%%%%%%%% 
\subsection{Second Order Markov Transition intensities}
%%%%%%%%%%%%%%%%%%%%%%%%%%%%%%%%%%%%%%%%%%%%%%%%%%%%%%%% 

Along with the definition of the second order transition probabilities we can also define the   second order transition intensities. 
\begin{definition}
For times $(s,t)\; s< t$ and adjacent states $h, j, k$, the  second order transition intensities are defined as 
\begin{eqnarray*}
\alpha_{hjk}(s,t)&=&\lim_{\Delta t\rightarrow 0} \frac{1}{\Delta t} P_{hjk}(s,t,t+\Delta t)
\end{eqnarray*}
That is, $\alpha_{hjk}(s,t)$ represents the propensity to arrive to the state $k$ at time $t$ knowing that in the previous instant the process was at state $j$ and that at time $s$ the process  was at state $h$. 
\end{definition}

As we have done for the transition probabilities we can express the transition intensities with $M$ matrices of dimension $M$, that is,  $\mathbf{\alpha_{(h)}}(s,t)=\left(\alpha_{hjk}(s,t)\right)_{j,k\in\mathcal{S}}$ for each $h\in\mathcal{S}$ and for each $(s,t),\,s<t$.
Under the homogeneous assumption we have:
\begin{eqnarray*}
\alpha_{hjk}(s,t)&=&\lim_{\Delta t\rightarrow 0} \frac{1}{\Delta t} P_{hjk}(t-s,\Delta t).
\end{eqnarray*}

%%%%%%%%%%%%%%%%%%%%%%%%%%%%%%%%%%%%%%%%%%%%%%%%%%%
 %%%%%%%%%%%%%%%%%%%%%%%%%%%%%%%%%%%%%%%%%%%%%%%%%%
\section{Extended Chapman-Kolmogorov equations} 
%%%%%%%%%%%%%%%%%%%%%%%%%%%%%%%%%%%%%%%%%%%%%%%%%%%
 %%%%%%%%%%%%%%%%%%%%%%%%%%%%%%%%%%%%%%%%%%%%%%%%%%
Given that clinical outcomes are often collected in days and aiming to compute the probability of being in a given state after a certain number of days, we
consider in  this section a discrete-time multistate process instead of a continuous-time stochastic process defined for $t\in[0,T]$.  
Other instances of discrete-time multistate process have been used to model
  COVID-19 disease progression and clinical outcomes (\cite{Chakladar}). We start reviewing Chapman-Kolmogorov equations for first order Markov chains to subsequently extend them to second order Markov chains.

%%%%%%%%%%%%%%%%%%%%%%%%%%%%%%%%%%%%%%%%%%%%%%%%%%%
\subsection{Chapman-Kolmogorov equations for first order Markov chains}
%%%%%%%%%%%%%%%%%%%%%%%%%%%%%%%%%%%%%%%%%%%%%%%%%%%
A first order discrete-time multistate models, known as Markov  chain, taking values in a discrete state space ${\mathcal S}=\{1,\cdots, M\}$ is the discrete version of a  
first order continuous-time Markov process. Hence, a
 Markov chain  is a stochastic model describing a sequence of possible events happening on discrete times in which the probability of each event depends only on the state attained in the previous event. 
 The Chapman-Kolmogorov relation is an important result in the theory of (\textbf{discrete}) Markov chains as it provides a method for calculating the  $n$-step transition probabilities.
The Chapman-Kolmogorov equations, for any    $s,t,u\in\mathbb{N} \; (s<u<t)$ are given by:
 \begin{equation}
      P_{hj}(s,t)=P_{hj}(s,u)\,P_{hj}(u,t)\label{eq:C-K_M}
 \end{equation}
and  follow as  a consequence of the Markov condition.
Chapman-Kolmogorov equations allow to reduce the general computation of  $P_{hj}(s,t)$, for any $s<t$,  ($s,t\in\mathbb{N})$ to the computation of one-time step transition probabilities,
$P_{hj}(s,s+1)$.
Denote by $\mathbf{P}(s)$ the one-time step transition probability matrix under the Markov assumption, that is,
  \begin{equation*}
    \mathbf{P}(s)=\{P_{hj}(s);\, h,j\in {\cal S}=\{1,\cdots, M\}\}
  \end{equation*}
 where $P_{hj}(s)$  stands for $P_{hj}(s,s+1)$. The collection of matrices $\mathbf{P}(s)$ is reduced to 
 the transition probability  matrix $  \mathbf{P}$ given by
  \begin{equation*}
    \mathbf{P}=\{P_{hj}=P_{hj}(1);\, h,j\in {\cal S}=\{1,\cdots, M\}\}
  \end{equation*}
 under the homogeneity  assumption (see Definition \ref{def:hom1}).
Hence, to study the evolution of the process for more than one time step, and thanks to the Chapman-Kolmogorov equations \eqref{eq:C-K_M}, it is only necessary to calculate the one-time initial transition probabilities. In the next section  we develop an extension of this result for second order Markov chains.

%%%%%%%%%%%%%%%%%%%%%%%%%%%%%%%%%%%%%%%%%%%%%%%%%%%
\subsection{Chapman-Kolmogorov equations for second order Markov chains} \label{Sec:CK2}
%%%%%%%%%%%%%%%%%%%%%%%%%%%%%%%%%%%%%%%%%%%%%%%%%%%
  
 Under second order Markov and homogeneity assumptions the transition probability matrices defined in Section \ref{sec:2norderTransProb} 
 satisfy, for times $s< t< u$ and {adjacent states} $h, j, k$:
 \begin{equation*}
P_{hjk}(s,t,u;{\mathcal{F}}_{s-})=P(X_u=k \mid X_s=h, X_t=j)=P_{hjk}(s,t,u) =P_{hjk}(t-s, u-t)
\end{equation*}
 
\noindent In particular, for  any $s\in\mathbb{N},\,s>1$ and consecutive times    $s, s+1, s+2$ we have 
 \begin{eqnarray*}
P_{hjk}(s,s+1,s+2;{\mathcal{F}}_{s-})&=&P(X_{s+2}=k \mid X_s=h, X_{s+1}=j)\\&=&P_{hjk}(s,s+1,s+2) =P_{hjk}(1, 1).
\end{eqnarray*}

\noindent Next theorem presents the equations to compute $n$-step transition probabilities such as
\begin{equation} \label{eq:nprob}
P(X_{s+n+1}=l|X_{s+1}=j,X_s=h)=P_{hjl}(1,n),
\end{equation}
for $ \;s,\,n\in\mathbb{N},\,n>1\;\mbox{and}\; l,j,k\in \mathcal{S}$ using only the  initial transition probabilities $P(X_3=l|X_2=j,X_1=h)=P_{hjl}(1,1)=P_{hjl}$.

Since 
 our desired probabilities  only depend on $n$, without lost of generality we can assume $s=1$ and the probabilities at \eqref{eq:nprob} can be equivalently written as
$$P(X_{n+2}=l|X_{2}=j,X_1=h),\quad \mbox{for} \quad n\in\mathbb{N},\; l,j,k\in \mathcal{S}
.$$

\noindent Recall that these transition probabilities can be summarized into $M$ matrices of dimension $M$ for each $n\in \mathbb{N}$. So for each state $h\in{\cal S}$ and $n>1$
\[\mathbf{P_{(h)}}(1,2,n+2)=\left(P_{hjl}(1,n)\right)_{j,k\in\mathcal{S}}=\left(P_{hjl}(n)\right)_{j,l\in\mathcal{S}},\]
we will omit $n$ of the previous notation when $n=1$.

\begin{notation} \label{not:mat}
Previous to the main result, we will present the matricial notation used in order to simplify the reading.

\begin{itemize}
\setlength{\itemsep}{-5pt}%
    \item Row $j$ of matrix $h$: $P_{j\cdot(h)}=(P_{hjk})_{k\in\mathcal{S}}$.
    \item Column $k$ of matrix $h$: $P_{\cdot k(h)}=(P_{hjk})_{j\in\mathcal{S}}$.
    \item We will use the $\ast$ symbol to denote that the elements of the row multiply each row of the matrix. For example $P_{hj\cdot}\ast P_{(h)}$ means that element $P_{hjk}$ multiplies all the elements of the row $k$ of $ P_{(h)}$.
    \item $\mathbf{P}^{(l)}$ is the matrix composed with the $l$ column of each of the $P_{h}$ matrices $h\in\mathcal{S}$.
    \item $\Tr()$ will denote the trace of a matrix.
\end{itemize}
\end{notation}

\begin{theorem}\label{th:C-K}
Assume $(X_n)_{n\in\mathbb{N}}$ is an homogeneous second order Markov chain. For any states $h,j,l\in \mathcal{S}$ and the notation defined in Notation \ref{not:mat}, we have
\begin{eqnarray}
P(X_{4}=l\mid X_{2}=j,X_1=h)&=&\mathbf{P}_{j\cdot (h)} \cdot \mathbf{P}_{\cdot l(j)} \label{eq:case1}\\
P(X_{5}=l\mid X_{2}=j,X_1=h)&=& \Tr\left({P_{hj\cdot}}\ast \mathbf{P}_{(j)}\cdot \mathbf{P}^{(l)}\right) \label{eq:case2}\\
P(X_{6}=l\mid X_{2}=j,X_1=h) &=& \sum_{k_1=1}^M P_{hjk_1}\cdot \Tr\left({P_{jk_1\cdot}}\ast \mathbf{P}_{(k_1)}\cdot \mathbf{P}^{(l)}\right) \label{eq:case3}\\
P(X_{7}=l\mid X_{2}=j,X_1=h) &=& \sum_{k_2=1}^M\sum_{k_1=1}^M P_{hjk_2}P_{jk_2k_1}\cdot \Tr\left({P_{k_2k_1\cdot}}\ast \mathbf{P}_{(k_1)}\cdot \mathbf{P}^{(l)}\right) \label{eq:case4}
\end{eqnarray}

General case $n\geq 7$
\begin{eqnarray*}
P(X_{n+1}=l\mid X_{2}=j,X_1=h) &=& \sum_{k_{n-4}=1}^M\ldots\sum_{k_2=1}^M\sum_{k_1=1}^M P_{hjk_{n-4}}P_{jk_{n-4}k_{n-3}}\ldots P_{k_3k_2k_1}\\ && \hspace{3cm}\times \Tr\left({P_{k_2k_1\cdot}}\ast \mathbf{P}_{(k_1)}\cdot \mathbf{P}^{(l)}\right)    
\end{eqnarray*}
\end{theorem}

\begin{proof}
The proof will be divided into three steps.

\textbf{Step 1.} We are proving the first case \eqref{eq:case1}.

Using the total probabilities Theorem and the second order Markov property,
\begin{eqnarray}\label{eq:calcase1}
&&\hspace{-1cm}P(X_{4}=l\mid  X_{2}=j,X_1=h)= \sum_{k=1}^M P(X_{4}=l, X_{3}=k\mid X_{2}=j,X_1=h)\nonumber\\
&& \hspace{-0.45cm}=\sum_{k=1}^M P(X_{3}=k \mid X_{2}=j,X_1=h) P(X_{4}=l \mid X_{3}=k,X_{2}=j,{X_1=h})\nonumber\\
&& \hspace{-0.45cm}=\sum_{k=1}^M P(X_{3}=k \mid X_{2}=j,X_1=h) P(X_{4}=l \mid X_{3}=k,X_{2}=j)\nonumber\\
&& \hspace{-0.45cm}=\sum_{k=1}^M  P_{hjk} \cdot {P_{jkl}}
\end{eqnarray}

We recall that $P(X_{4}=l \mid X_{3}=k,X_{2}=j)=P_{jkl}(1,1)=P_{jkl}$ since for all the assumptions only the difference between times determine the transition probabilities. From here we can write it in matricial form as in \eqref{eq:case1}.

\textbf{Step 2.} Now we focus in the second case \eqref{eq:case2}.

Using the same argument of the previous case and also the result obtained \eqref{eq:calcase1}
\begin{eqnarray*}
&&\hspace{-1cm}P(X_{5}=l\mid  X_{2}=j,X_1=h)= \sum_{k_1=1}^M P(X_{5}=l, X_{3}=k_1\mid  X_{2}=j,X_1=h)\nonumber\\
&& \hspace{-0.45cm}=\sum_{k_1=1}^M P(X_{5}=l\mid X_{3}=k_1,X_{2}=j,{X_1=i})P(X_{3}=k_1\mid  X_{2}=j,X_1=h)\nonumber\\
&& \hspace{-0.45cm}=\sum_{k_1=1}^M P(X_{5}=l\mid X_{3}=k_1,X_{2}=j)P(X_{3}=k_1\mid  X_{2}=j,X_1=h)\nonumber\\
&& \hspace{-0.45cm}=\sum_{k_1=1}^M \left(\sum_{k_2=1}^M{P_{k_1k_2l}}\cdot {P_{jk_1k_2} }\right)\cdot {P_{hjk_1}}
\end{eqnarray*}

Now, if we want to write it in a matricial way we can observe that $\sum_{k_2=1}^M{P_{k_1k_2l}}\cdot {P_{jk_1k_2} }$ corresponds to the product of the $l$ column of each matrix by the matrix $\mathbf{P}_{(j)}$ and then each row of this matrix product is multiplied by the probabilities ${P_{hjk_1}}$ which are the elements of row $j$ of matrix $\mathbf{P}_{(h)}$. From here we obtain the formula in \eqref{eq:case2}.

\textbf{Step 3.} We follow by proving \eqref{eq:case3}.

We repeat here the arguments in the previous steps and also we apply the previous results.
\begin{eqnarray*}
&&\hspace{-1cm}P(X_{6}=l\mid  X_{2}=j,X_1=h)= \sum_{k_3=1}^M P(X_{6}=l, X_{3}=k_3\mid  X_{2}=j,X_1=h)\nonumber\\
&& \hspace{-0.45cm}=\sum_{k_3=1}^M P(X_{6}=l\mid X_{3}=k_3,X_{2}=j,{X_1=h})P(X_{3}=k_3\mid  X_{2}=j,X_1=h)\nonumber\\
&& \hspace{-0.45cm}=\sum_{k_3=1}^M P(X_{6}=l\mid X_{3}=k_3,X_{2}=j)P(X_{3}=k_3\mid  X_{2}=j,X_1=h)\nonumber\\
&& \hspace{-0.45cm}=\sum_{k_3=1}^M P_{hjk_3} \left(\sum_{k_1=1}^M{P_{jk_3k_1}} \sum_{k_2=1}^M{P_{k_1k_2l}}\cdot {P_{k_3k_1k_2}} \right)\\
&& \hspace{-0.45cm}= P_{hj1} \left(\sum_{k_1=1}^M{P_{jk_3k_1}} \sum_{k_2=1}^M{P_{k_1k_2l}}\cdot {P_{k_3k_1k_2}} \right)+\ldots+P_{hjM} \left(\sum_{k_1=1}^M{P_{jk_3k_1}} \sum_{k_2=1}^M{P_{k_1k_2l}}\cdot {P_{k_3k_1k_2}} \right)
\end{eqnarray*}

We observe that we obtain a similar expression as the previous step but multiplied by the probabilities $P_{hj\cdot}$ that correspond to the row $j$ of matrix $\mathbf{P}_{(h)}$.Thus, the matricial expression for this case is immediate.

From here, in order to prove Equation \eqref{eq:case4} and the general case we can just repeat the same arguments as in this last case to easily obtain the general formula by induction. In these two cases the principal modifications of the matricial form will focus in adding one sum for each step.
\end{proof}

\begin{corollary}\label{corChK}
Assume $(X_n)_{n\in\mathbb{N}}$ is an homogeneous second order Markov chain. For time $s>0$, any states $h,j,l\in \mathcal{S}$ and the notation defined in Notation \ref{not:mat}, we have
\begin{eqnarray*}
P(X_{s+3}=l\mid X_{s+1}=j,X_s=h)&=&\mathbf{P}_{j\cdot (h)} \cdot \mathbf{P}_{\cdot l(j)} \\
P(X_{s+4}=l\mid X_{s+1}=j,X_s=h)&=& \Tr\left({P_{hj\cdot}}\ast \mathbf{P}_{(j)}\cdot \mathbf{P}^{(l)}\right) \\
P(X_{s+5}=l\mid X_{s+1}=j,X_s=h) &=& \sum_{k_1=1}^M P_{hjk_1}\cdot \Tr\left({P_{jk_1\cdot}}\ast \mathbf{P}_{(k_1)}\cdot \mathbf{P}^{(l)}\right) \\
P(X_{s+6}=l\mid X_{s+1}=j,X_s=h) &=& \sum_{k_2=1}^M\sum_{k_1=1}^M P_{hjk_2}P_{jk_2k_1}\cdot \Tr\left({P_{k_2k_1\cdot}}\ast \mathbf{P}_{(k_1)}\cdot \mathbf{P}^{(l)}\right) 
\end{eqnarray*}

General case $n\geq 7$
\begin{eqnarray*}
P(X_{s+n}=l\mid X_{s+1}=j,X_s=h) &=& \sum_{k_{n-4}=1}^M\ldots\sum_{k_2=1}^M\sum_{k_1=1}^M P_{hjk_{n-4}}P_{jk_{n-4}k_{n-3}}\ldots P_{k_3k_2k_1}\\ && \hspace{3cm}\times \Tr\left({P_{k_2k_1\cdot}}\ast \mathbf{P}_{(k_1)}\cdot \mathbf{P}^{(l)}\right)    
\end{eqnarray*}
\end{corollary}

\begin{remark}
We observe that the extended Chapman-Kolmogorov equations only consider the case where the two past times are consecutive. If they are not consecutive, since we cannot know which is the intermediate state, we cannot provide any formula to compute it. For some models and specific cases with non return states it is possible to study it.
\end{remark}
%%%%%%%%%%%%%%%%%%%%%%%%%%%%%%%%%%%%%%%%%%%%%%%%%%%
%%%%%%%%%%%%%%%%%%%%%%%%%%%%%%%%%%%%%%%%%%%%%%%%%%%
\section{Estimation and inference under second order Markov assumption} \label{Sec:Estimation}
%%%%%%%%%%%%%%%%%%%%%%%%%%%%%%%%%%%%%%%%%%%%%%%%%%%
%%%%%%%%%%%%%%%%%%%%%%%%%%%%%%%%%%%%%%%%%%%%%%%%%%%
Given three different adjacent states $h,\,j,\,l\in\mathcal{S}$ ($h\rightarrow j \rightarrow l$) such that $h,\,j$ are not absorbent,
the purpose of this Section is to estimate
the $r$-step transition probabilities $P(X_{s+r}=l\mid X_{s-1}=j,X_{s-2}=h)$ for  
any $s,r\in\mathbb{N},\,s,r>1$.  
Under the homogeneity assumption we have  that $P(X_{s}=l\mid X_{s-1}=j,X_{s-2}=h)=P_{hjl}(1,1)=P_{hjl}$ and,
following Corollary \ref{corChK},  in order to estimate $P(X_{s+r}=l\mid X_{s-1}=j,X_{s-2}=h)$ is enough to estimate  the initial transition probabilities $P_{hjl}(1,1)=P_{hjl}$.

We assume that individuals  are followed  until  a maximum of $T$ units of time (days as in the illustration). Let $\{X^i_s,\ s=0, 1,\cdots, T\}$ denote the non-reversible multistate process for subject $i=1,\,\ldots,\,n$, where $X^i_s\in\mathcal{S}$. 
For $i=1,\,\ldots,\,n$, $s=2,\cdots, T$ and $h,\, j,\,l\in\mathcal{S}$  we define the counting processes \[N^i_{hjl}(s)=\mathbb{1}{\{X^i_{s-2}=h,\, X^i_{s-1}=j,\, X^i_s=l\}}\]
counting  $1$ if  subject $i$  has transit from state $h$ to state $j$ and to state $l$ at times $s-2,\,s-1,\,s$, respectively; and $0$ otherwise.  The total number of individuals who have  followed the path  $h\rightarrow j \rightarrow l$ at times $s-2,\,s-1,\,s$ is given by the sum $\widetilde N_{hjl}(s)= \sum_{i=1}^n N^i_{hjl}(s)$.  $\widetilde N_{hjl}(s)$
 is a binomial random variable with parameters $(n, \pi_{hjl}(s))$ where the  probability $\pi_{hjl}(s)$ corresponds to 
$\pi_{hjl}(s)=P(X_{s-2}=h,\, X_{s-1}=j,\, X_s=l)$. 
We also define the at-risk process of subject $i$ corresponding to states $h$ and $j$ at times $s=2, \cdots, T$, 
\[Y^i_{hj}(s-1)=\mathbb{1}{\{X^i_{s-2}=h,\, X^i_{s-1}=j\}}\]
counting $1$ if subject $i$ was at risk of moving to adjacent states to $j$ given that he/she was in states $h$ and $j$ at times $s-2$ and $s-1$,  respectively. The total number of individuals at risk at time $s$ is given by
 $\widetilde Y_{hj}(s-1)= \sum_{i=1}^n Y^i_{hj}(s-1)$ and corresponds to a  binomial random variable with parameters $(n, \pi_{hj}(s-1))$ where 
 $\pi_{hj}(s-1)=P(X_{s-2}=h,\, X_{s-1}=j)$.

Regarding the  estimation of the transition probability $P_{hjl}(s)=P(X_{s}=l\mid X_{s-1}=j, X_{s-2}=h)$ for a given $s\geq 2$,  we will proceed in two different ways. The first one  takes advantage of the ratio of the two probabilities
$P_{hjl}(s)={\pi_{hjl}(s)}/{\pi_{hj}(s-1)}$ while the second one exploits directly the estimation of the conditional probability $P_{hjl}(s)$.

For given adjacent states $(h, j)$ we have to take into account for which times $s$ there are individuals at risk. For this reason we introduce here the notation  $R_{hj}=\inf\{s\geq 2: \widetilde Y_{hj}(s-1)>0\}$ and
 $T_{hj}=\sup\{s\geq 2: \widetilde Y_{hj}(s-1)>0\}$ indicating that for any $s$ such that $R_{hj}\leq s\leq T_{hj}$, the estimation of $P_{hjl}$, for $l$ adjacent to $j$, would be plausible.  
 From a practical point of view we will have to guarantee that the number of individuals at risk $\widetilde Y_{hj}(s-1)$ is large enough for a meaningful estimation of $P_{hjl}$. It would be equivalent if we consider the number of days within each couple of states. In this case $R'_{hj}=2$ and  $T'_{hj}=T_{hj}-R_{hj}$.
%%%%%%%%%%%%%%%%%%%%%%%%%%%%%%%%%%%%%%%%%%%%%%%%%%%%
 \subsection{Estimation of $P_{hjl}$ via the Bernoulli probabilities $\pi_{hjl}$ and $\pi_{hj}$}
 %%%%%%%%%%%%%%%%%%%%%%%%%%%%%%%%%%%%%%%%%%%%%%%%%%

Given that  for all $s=2, \cdots, T$
\begin{eqnarray*} 
P_{hjl}&=&P_{hjl}(s)=P(X_{s}=l\mid X_{s-1}=j, X_{s-2}=h)=\frac{P(X_{s-2}=h,  X_{s-1}=j, X_{s}=l)}{ P(X_{s-2}=h, X_{s-1}=j)}\\
&=&\frac{\pi_{hjl}(s)}{\pi_{hj}(s-1)}
\end{eqnarray*}
a natural estimator  for $P_{hjl}$ can be built 
estimating separately both numerator and denominator by
$\sum_{s=R_{hj}}^{T_{hj}} \widetilde N_{hjl}(s)/n$ and $\sum_{s=R_{hj}}^{T_{hj}}\widetilde Y_{hj}(s-1)/n$, respectively.
Observe that $\sum_{s=R_{hj}}^{T_{hj}} \widetilde N_{hjl}(s)$
corresponds to the total number of individuals that have  followed the path  $h\rightarrow j \rightarrow l$ at any three times $(s-2,s-1,s)$ and $\sum_{s=R_{hj}}^{T_{hj}}\widetilde Y_{hj}(s-1)$
is the total number of individuals that have followed the path $h\rightarrow j$ consecutively at any two times $(s-2,s-1)$. 
\begin{theorem}
For given adjacent states $(h,j,l)$, the statistic
\begin{equation}
\widetilde{P}_{hjl} =\dfrac{\sum_{s=R_{hj}}^{T_{hj}} \widetilde N_{hjl}(s)}{\sum_{s=R_{hj}}^{T_{hj}}\widetilde Y_{hj}(s-1)}
\end{equation}
estimates consistently $P_{hjl}$.
\end{theorem}

\begin{proof}
Note that by the Law of Large Numbers we have convergence in probability of the following two  estimators:
\begin{eqnarray*} 
&&\frac{1}{n}\sum_{s=R_{hj}}^{T_{hj}} \widetilde N_{hjl}(s)=\frac{1}{n}\sum_{i=1}^n\sum_{s=R_{hj}}^{T_{hj}} N^i_{hjl}(s)\xrightarrow[n\rightarrow\infty]{\text{P}} E\left(\sum_{s=R_{hj}}^{T_{hj}} N^i_{hjl}(s)\right)=\sum_{s=R_{hj}}^{T_{hj}}E\left( N^i_{jhl}(s)\right)\\&&\hspace{9cm}=\sum_{s=R_{hj}}^{T_{hj}}\pi_{hjl}(s)\\
&&\frac{1}{n}\sum_{s=R_{hj}}^{T_{hj}} \widetilde Y_{hj}(s-1)=\frac{1}{n}\sum_{i=1}^n\sum_{s=R_{hj}}^{T_{hj}} Y^i_{hj}(s-1)\xrightarrow[n\rightarrow\infty]{\text{P}}E\left(\sum_{s=R_{hj}}^{T_{hj}} Y^i_{hj}(s-1)\right)\\&&\hspace{6.5cm}=\sum_{s=R_{hj}}^{T_{hj}}E\left( Y^i_{jh}(s-1)\right)=\sum_{s=R_{hj}}^{T_{hj}}\pi_{hj}(s-1)
\end{eqnarray*}

Second order homogeneity implies that
\begin{eqnarray*}
    \pi_{hjl}(s)&=& P\{X_{s-2}=h,X_{s-1}=j, X_s=l\}\\
    &=&P\{X_s=l \mid X_{s-2}=h, X_{s-1}=j\}P\{X_{s-1}=j,  X_{s-2}=h\}
    \\
   &= &P_{hjl}(1,1)P\{X_{s-1}=j, X_{s-2}=h\}= P_{hjl}(1,1)\pi_{hj}(s-1),
\end{eqnarray*}
 and we conclude that $\widetilde{P_{jhl}}$ converges in probability to $P_{hjl}$:
\begin{eqnarray*}
\widetilde{P}_{jhl} & =&\dfrac{\sum_{s=R_{hj}}^{T_{hj}} \tilde N_{hjl}(s)}{\sum_{s=R_{hj}}^{T_{hj}}\tilde Y_{hj}(s-1)}=
\dfrac{\sum_{s=R_{hj}}^{T_{hj}} \tilde N_{hjl}(s)/n}{\sum_{s=R_{hj}}^{T_{hj}}\tilde Y_{hj}(s-1)/n} \xrightarrow[n\rightarrow\infty]{\text{P}} \dfrac{\sum_{s=R_{hj}}^{T_{hj}}\pi_{hjl}(s)}{\sum_{s=R_{hj}}^{T_{hj}}\pi_{hj}(s-1)}
\\&=&\dfrac{\sum_{s=R_{hj}}^{T_{hj}}P_{hjl}(1,1)\pi_{hj}(s-1)}{\sum_{s=R_{hj}}^{T_{hj}}\pi_{hj}(s-1)}=
P_{hjl}
\end{eqnarray*}
\end{proof}

 %%%%%%%%%%%%%%%%%%%%%%%%%%%%%%%%%%%%%%%%%%%%%%%%%%%%%%%%%%
\subsection{Estimation of $P_{hjl}$ via the conditional probability}\label{sec:EstCondProb}
%%%%%%%%%%%%%%%%%%%%%%%%%%%%%%%%%%%%%%%%%%%%%%%%%%%%%%%%%%
 For every $s$ such that $R_{hj}\leq s\leq T_{hj}$,
the relative frequency given by the ratio $$\widehat{P}_{hjl}(s)=\dfrac{\ \widetilde N_{hjl}(s)}{ \widetilde Y_{hj}(s-1)}$$ is an straightforward estimator of $P_{hjl}(s)$.
Since for all $s=2, \cdots, T$, $P_{hjl}(s)=P_{hjl}$  an estimator for $P_{hjl}$ can be obtained as the average of $\widehat{P}_{hjl}(s)$, that is,
\begin{equation}\label{def:Phat}
\widehat{P}_{hjl} =\frac{1}{T_{hj}-R_{hj}+1}\sum_{s=R_{hj}}^{T_{hj}}
\widehat{P}_{hjl}(s)=\frac{1}{T_{hj}-R_{hj}+1}\sum_{s=R_{hj}}^{T_{hj}}
\dfrac{\widetilde N_{hjl}(s)}{ \widetilde Y_{hj}(s-1)}. 
\end{equation}

\begin{theorem}
For given adjacent states $(h,j,l)$ the statistic
$\widehat{P}_{hjl}$ given in \eqref{def:Phat} is an unbiased estimatator of $P_{hjl}$.
\end{theorem}
\begin{proof}
    
We can write $\widehat{P}_{hjl}(s)$  as \begin{equation*}
\widehat{P}_{hjl}(s)=\dfrac{\ \widetilde N_{hjl}(s)}{ \widetilde Y_{hj}(s-1)} =
\sum_{\{i: X^i_{s-2}=h, X^i_{s-1}=j\}}\dfrac{\ N^i_{hjl}(s)}{ \widetilde Y_{hj}(s-1)}. 
\end{equation*}
$N^i_{hjl}(s)$, conditioned to those individuals  being consecutively in states $h, j$  at times $s-2,$ $s-1$, is a Bernoulli random variable with probability $P_{hjl}(s)$ and $\widetilde N_{hjl}(s)$ is a Binomial random variable with number of trials equal to $\widetilde Y_{hj}(s-1) $ and probability of success equal to $P_{hjl}(s)=P_{hjl}$. Hence, 
$E(\widehat{P}_{hjl}(s)|\widetilde Y_{hj}(s-1))=P_{hjl}$. Note that for all s, 
\begin{equation*}
E\left(\widehat{P}_{hjl}(s)\right)=E\left\{E\left(\widehat{P}_{hjl}(s)|\widetilde Y_{hj}(s-1)\right)\right\}=E\left\{{P}_{hjl}(s)\right\}=P_{hjl}
\end{equation*}
Somehow heuristically, considering $R_{hj}$ and $T_{hj}$ fixed not random values, 
\begin{equation*}
E(\widehat{P}_{hjl} )=E\left(\frac{1}{T_{hj}-R_{hj}+1}\sum_{s=R_{hj}}^{T_{hj}}
\widehat{P}_{hjl}(s)\right)=\frac{1}{T_{hj}-R_{hj}+1}\sum_{s=R_{hj}}^{T_{hj}}
E\left(\widehat{P}_{hjl}(s)\right)=P_{hjl}.
\end{equation*}
\end{proof}
%%%%%%%%%%%%%%%%%%%%%%%%%%%%%%%%%%%%%
%%%%%%%%%%%%%%%%%%%%%%%%%%%%%%%%%%%%%
\section{DIVINE model}
%%%%%%%%%%%%%%%%%%%%%%%%%%%%%%%%%%%%%
%%%%%%%%%%%%%%%%%%%%%%%%%%%%%%%%%%%%%

%%%%%%%%%%%%%%%%%%%%%%%%%%%%%%%%%%%%%
\subsection{Description}
%%%%%%%%%%%%%%%%%%%%%%%%%%%%%%%%%%%%%
The dataset we use  as illustration corresponds to  a cohort of 2076 COVID-19 hospitalised patients (during the first wave of the pandemic, March-April 2020) in five hospitals located in the southern Barcelona metropolitan area (Spain). 
Since all the patients were  monitored until discharge from hospital or death, the transition times (in days) are known exactly for all subjects and there are not incomplete data due to lost to follow-up.

\begin{figure}[h!] 
\centering
\includegraphics[width=0.9\textwidth]{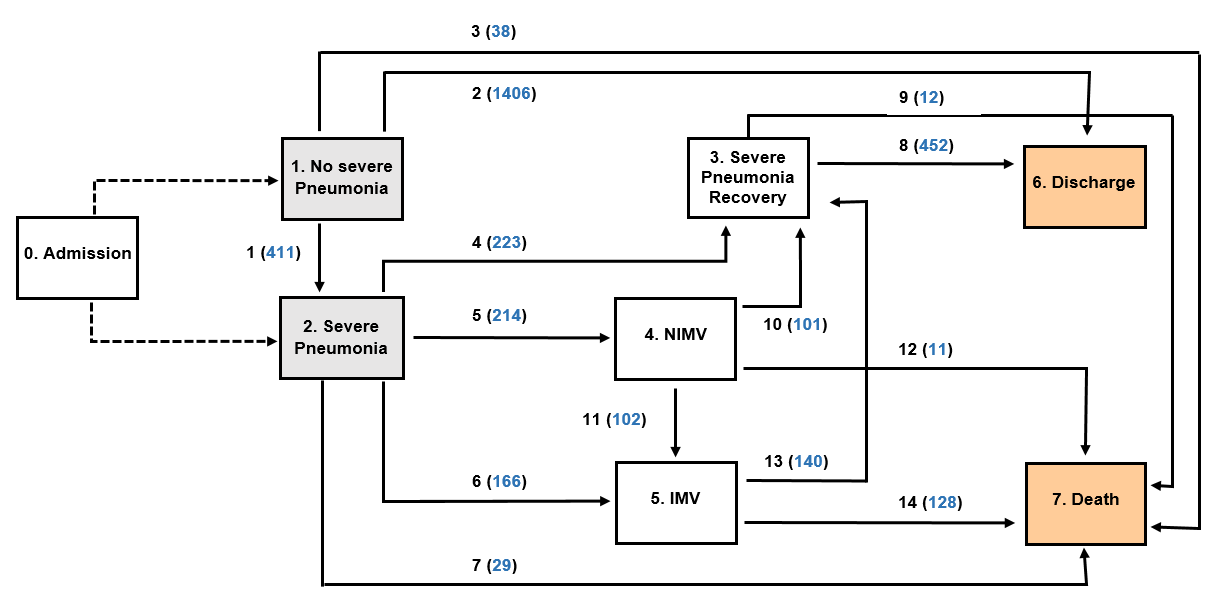}
\caption{Graphical representation of the multistate model for modelling the trajectory of hospitalized COVID-19 patients. Seven states are considered and 14 possible transitions (in parentheses the sample size of each transition).
\footnotesize{NIMV: Non-Invasive mechanical ventilation, IMV: Invasive mechanical ventilation.}}
\label{FullModel}
\end{figure}

This data is part of the DIVINE project\\ (\url{https://grbio.upc.edu/en/research/highlighted-projects})\\ for which a  
 multidisciplinary research team integrated by researchers from the GRBIO (UPC-UB), Bellvitge University Hospital, and Bellvitge Biomedical Research Institute has collaborated to define a statistical framework with a clear clinician focus on achieving deeper understanding of the severe form of the disease caused by the SARS-CoV-2 virus. Based on the  team cooperative knowledge 
 a multistate model with seven states and 14 transitions has been built (see Figure \ref{FullModel} where the numbers in parentheses denote the patients doing that transition). As seen in Figure \ref{FullModel},  7  states are considered: (1) No Severe Pneumonia (NSP), (2) Severe Pneumonia (SP), (3) Severe Pneumonia Recovery (Recov), (4) Non invasive mechanical ventilation (NIMV), (5) Invasive mechanical ventilation (IMV), (6) Discharge (Disch) and (7) Death (Death).

The following considerations are in place: 
\begin{enumerate}
    \item Once a patient has been  admitted (state 0), he/she is immediately assigned to one of the  two initial states: No Severe Pneumonia and Severe Pneumonia. It is assumed that the process starts at time  $t=0$ in one of these  two states.
    \item Discharge and Death are absorbing states implying that once a patient has been   discharged or has died  he/she cannot re-enter to be hospitalised again.  
    \item The time scale used in this model is days since the hospital admission. For all the transitions, the transition times (in days) are exactly known.
    \item Patients can only jump to a neighboring state in a single day.
\end{enumerate}
For more details on the data and the clinical patient characteristics see \cite{Pallares, Garmendia, Piulachs}. 

The main goal with this illustration is to study the evolution of the patients without the restriction of a first order Markov assumption.  To do so we start 
 validating  for which transitions of the previous multistate model the Markov assumption holds. Next, we will estimate 
 the transition probabilities between two states taking into account that they might depend as well  on the  immediate previous state. Finally, we will  compare the evolution of those patients admitted with No Severe Pneumonia versus those admitted with Severe Pneumonia.

%%%%%%%%%%%%%%%%%%%%%%%%%%%%%%%%%%%%%%%%%%%%%%%
\subsection{Description of direct and two-step transitions}
%%%%%%%%%%%%%%%%%%%%%%%%%%%%%%%%%%%%%%%%%%%%%%%

Table \ref{tab:tab1} summarises the number of patients for each direct transition and the number of patients for the corresponding related 2-step transitions (consecutive states but not necessarily consecutive times).  For instance, individuals doing the direct transition  $\mbox{Recov}\rightarrow \mbox{Disch}$  might arrive from $\mbox{SP}:$ $\mbox{SP}\rightarrow \mbox{Recov} \rightarrow \mbox{Disch}$ from $\mbox{NIMV}:$ $ \mbox{NIMV}\rightarrow \mbox{Recov} \rightarrow \mbox{Disch}$  or
from $\mbox{IMV}:$ $ \mbox{IMV}\rightarrow \mbox{Recov}\rightarrow \mbox{Disch}$. Note that we are only considering those direct transitions $j\rightarrow l$ for which there exists, at least,  an state $k$ adjacent to $j$ ($k\rightarrow j\rightarrow l$).

\begin{table}[H]
\centering
    \begin{tabular}{|l|c|r|c|c|}
    \hline
    Direct transition & Sample size & 2-step transition & Sample size & Percent.\\
    \hline
        \multirow{2}{*}{$ \mbox{SP} \rightarrow \mbox{Recov}$}& \multirow{2}{*}{223} &$\mbox{NSP}\rightarrow \mbox{SP} \rightarrow \mbox{Recov}$ & 171  &76.68 \\
        &&$\mbox{SP} \rightarrow \mbox{Recov}$ & 52  &23.32\\
        \hline
        \multirow{2}{*}{$ \mbox{SP} \rightarrow \mbox{NIMV}$}& \multirow{2}{*}{214} &$\mbox{NSP}\rightarrow \mbox{SP} \rightarrow \mbox{NIMV}$ & 134  & 62.62  \\
        &&$\mbox{SP} \rightarrow \mbox{NIMV}$ & 80 & 37.38\\
        \hline
         \multirow{2}{*}{$ \mbox{SP} \rightarrow \mbox{IMV}$}& \multirow{2}{*}{166}&$\mbox{NSP}\rightarrow \mbox{SP} \rightarrow \mbox{IMV}$ & 92 & 55.42\\
        &&$\qquad \mbox{SP} \rightarrow \mbox{IMV}$ & 74 & 44.58\\
        \hline
        \multirow{2}{*}{$ \mbox{SP} \rightarrow \mbox{Death}$}& \multirow{2}{*}{29} &$\mbox{NSP}\rightarrow \mbox{SP} \rightarrow \mbox{Death}$ & 14 &48.26 \\
        &&$\qquad \mbox{SP} \rightarrow \mbox{Death}$ &  15 &51.72 \\
        \hline
        \multirow{3}{*}{$\mbox{Recov} \rightarrow \mbox{Disch}$}& \multirow{3}{*}{452} &$\mbox{SP}\rightarrow \mbox{Recov} \rightarrow \mbox{Disch}$ & 223  &49.34 \\
        &&$\mbox{NIMV}\rightarrow \mbox{Recov} \rightarrow \mbox{Disch}$ & 96 &21.24\\
        &&$\mbox{IMV}\rightarrow \mbox{Recov} \rightarrow \mbox{Disch}$ & 133 &29.42\\
        \hline
        \multirow{3}{*}{$\mbox{Recov} \rightarrow \mbox{Death}$}& \multirow{3}{*}{12} &$\mbox{SP}\rightarrow \mbox{Recov} \rightarrow \mbox{Death}$ & 0 & 0\\
        &&$\mbox{NIMV}\rightarrow \mbox{Recov} \rightarrow \mbox{Death}$ & 5 &41.67\\
        &&$\mbox{IMV}\rightarrow \mbox{Recov} \rightarrow \mbox{Death}$ & 7 &58.33\\
        \hline
        \multirow{2}{*}{$\mbox{IMV} \rightarrow \mbox{Death}$}& \multirow{2}{*}{128}&$\mbox{SP}\rightarrow \mbox{IMV} \rightarrow \mbox{Death}$ & 71&57.26 \\
        &&$\mbox{NIMV}\rightarrow \mbox{IMV} \rightarrow \mbox{Death}$ & 57 &45.97\\
        \hline
        \multirow{2}{*}{$\mbox{IMV} \rightarrow \mbox{Recov}$}& \multirow{2}{*}{140} &$\mbox{SP}\rightarrow \mbox{IMV} \rightarrow \mbox{Recov}$ & 95 &67.86\\
        &&$\mbox{NIMV}\rightarrow \mbox{IMV} \rightarrow \mbox{Recov}$ & 45 &32.14\\\hline
    \end{tabular}\\
    \parbox{13.5cm}{\footnotesize{States: NSP: No Severe Pneunomia, SP: Severe pneumonia, Recov: Severe Pneumonia recovery, NIMV: Non-Invasive Mechanical Ventilation, IMV: Invasive Mechanical Ventilation, Disch: Discharge, Death: Death}}
    \caption{Aggregation of the 2 steps paths for each direct transition taking into account the previous immediate state} 
    \label{tab:tab1}
\end{table}

Relating the sample size of the direct transitions appearing in Table \ref{tab:tab1} with the notation introduced in Section \ref{Sec:Estimation}, we see that the sample size of the direct transition $h\rightarrow j$  corresponds to the  number of patients at-risk: $\sum_{s=2}^T \tilde{Y}_{hj}$, while  the sample size of the 2-step transition $h\rightarrow j\rightarrow l$ does not coincide with $\sum_{s=2}^T \tilde{N}_{hjl}$  since we have considered all patients doing this path either in consecutive times or not.

Table \ref{tab:tab1} reveals that  the proportion of patients for a given  transition (e.g,   IMV$\rightarrow$ Recov)  drastically differs whether the patients were before in SP ($68\%$) or in NIMV ($32\%$). We  can also examine the transition SP$\rightarrow$Recov, if we separate the patients between those with NSP in the admission (76.68\%) and those with SP in the admission (23.32\%) we also observe important differences. Similar interpretation is in place with transition SP$\rightarrow$NIMV. This suggests that the model may not fulfill the Markov assumption and that it may be important to take the two previous states into account when calculating the transition probabilities.

%%%%%%%%%%%%%%%%%%%%%%%%%%%%%%%%%%%%%%%%%%%%%%%
\subsection{Testing the Markov assumption} \label{Sec: MarkovTest2}
%%%%%%%%%%%%%%%%%%%%%%%%%%%%%%%%%%%%%%%%%%%%%%%
In order to check which 2-step transitions are not first order Markovian, we use the Markov test described in Section \ref{sec:MarkovTest} and follow  \cite{TitmanPutter} guidelines with respect to the 
time intervals  $[t_0,t_{\max}]$  where the test can be conducted. Basically,  the comparison is restricted to windows of time with enough individuals and to direct transitions that have an immediate previous state.

To evaluate the logrank test we compute the statistics for an equally 0.5-day spaced grid in the interval $[1,11]$ for all the transitions except for transitions 7 (SP $\rightarrow$ Death) and 12 (NIMV $\rightarrow$ Death) in which the interval is $[1,7]$ and transition 14 (IMV $\rightarrow$ Death) with the interval $[1,16]$.

Table \ref{tab:MarkovTest} summarizes the p-values of the log-rank tests obtained from 5000 wild bootstrap resamples (\cite{Lin(1993)}) and considering the three possible summary statistics: weighted mean, mean and supremum described in Section \ref{sec:MarkovTest}. For each transition (rows), we have carried out the test for all the possible previous states as well as the overall chi-squared test. 
The partial p-values are the ones corresponding to the global test. For each transition, the transition intensity compares  the subjects who were previously at fixed state $j$ (in columns) versus the ones who were not there.
\begin{table}[H]
    \centering
\begin{tabular}{l|ccccccc}
\hline
Transitions&& NSP & SP & Recov & NIMV & IMV &overall\\
\hline
4 \;\,(SP$\rightarrow$ Recov) & UM & 0.005& 0.005 &&&&\textbf{0.0042}\\
 & WM & 0.006 & 0.006 &&&&\\
  & S & 0.029 & 0.029 &&&&\\
5 \;\,(SP$\rightarrow$ NIMV) &UM &$<10^{-16}$ & $<10^{-16}$ &&&&\textbf{0.0018}\\
 &WM &$<10^{-16}$ & $<10^{-16}$ &&&&\\
 &S &0.026 & 0.026 &&&&\\
6 \;\,(SP$\rightarrow$ IMV) &UM & 0.009 & 0.009 &&&&\textbf{0.016}\\ 
&WM & 0.002 & 0.002 &&&&\\
&S & 0.042 & 0.042 &&&&\\
7 \;\,(SP$\rightarrow$ Death) &UM &0.106 &0.106 &&&&0.196\\
&WM &0.120 &0.120 &&&&\\
& S &0.340 &0.340 &&&&\\
8 \;\,(Recov$\rightarrow$ Disch) &UM &0.007 & $<10^{-5}$ &0.165 & $<10^{-5}$ & $<10^{-5}$ &$\mathbf{<10^{-16}}$\\
 &WM &0.010 & $<10^{-5}$ &0.186 & $<10^{-5}$ & $<10^{-5}$ &\\
  &S &0.104 & $<10^{-5}$ &0.388 & $<10^{-5}$ & $<10^{-5}$ &\\
9 \;\,(Recov$\rightarrow$ Death) &UM& 0.652&0.298 &0.145 &0.495 &0.143 &0.357\\
&WM& 0.644&0.273 &0.144 &0.464 &0.151 &\\
&S& 0.656&0.313 &0.353 &0.639 &0.309 &\\
10 (NIMV$\rightarrow$ Recov) &UM&0.594&0.190&&0.694&&0.609\\
&WM&0.588&0.183&&0.717&&\\
&S&0.831&0.432&&0.892&&\\
11 (NIMV$\rightarrow$IMV) &UM&0.514& 0.819& &0.728&&0.807\\
&WM&0.501& 0.858& &0.765&&\\
&S&0.432& 0.311& &0.304&&\\
12 (NIMV$\rightarrow$ Death) &UM &0.348&0.218&&0.649&&0.437\\
&WM &0.378&0.253&&0.619&&\\
&S &0.338&0.342&&0.719&&\\
13 (IMV$\rightarrow$ Recov) &UM &0.564&$<10^{-3}$& &0.514&0.005&0.059\\
&WM &0.531&$<10^{-3}$& &0.456&0.005&\\
&S &0.780&0.034& &0.376&0.037&\\
14 (IMV$\rightarrow$ Death) &UM &0.296&0.663&&0.318&0.099&0.305\\
&WM &0.296&0.674&&0.269&0.100&\\
&S &0.471&0.369&&0.264&0.205&\\
\hline
\end{tabular}
\caption{p-values obtained from the computation of the Markov test for each transition and each previous state. Three different summary statistics have been computed: unweighted mean (UM), weighted mean (WM) and supremum (S). In bold the transitions that are statistically significant at 0.05.} \label{tab:MarkovTest}
\end{table}

 Considering the overall p-values, transitions 4 (SP $\rightarrow$ Recov), 5 (SP $\rightarrow$ NIMV), 6 (SP $\rightarrow$ IMV) and  8 (Recov $\rightarrow$ Disch) show clear departures from the Markov assumption, while transition 13 (IMV $\rightarrow$ Recov) is marginally significant. Furthermore,  any one of the summary tests rejects the Markovianity in transitions 4, 5 and 6. The supremum statistic would not reject Markovianity  from states NSP and Recov in transition 8. Finally, the global p-value of 0.059 in transition 13 is mainly due to the non Markovianity coming from states SP and IMV.

 These findings suggest that once a patient is critically ill, for instance, in states NIMV and IMV, the future clinical evolution is independent of whether he/she was diagnosed with NSP or SP when hospitalized. However, the clinical evolution to NIMV or IMV will be different for those patients initially diagnosed with NSP versus those diagnosed with SP when hospitalized.
These results lead us to  consider second order Markov multistate models  in order to study the evolution of the hospitalized COVID-19 patients during the fist wave of the pandemia. 

%%%%%%%%%%%%%%%%%%%%%%%%%%%%%%%%%%%%%%%%%%%%%%%
\subsection{Estimation of the transition probability matrices}
%%%%%%%%%%%%%%%%%%%%%%%%%%%%%%%%%%%%%%%%%%%%%%%
We will now estimate the seven possible transition probability matrices using the estimators presented in Section \ref{Sec:Estimation}. As we have mentioned before, we know the exact transition times, so we can easily estimate the transition probability by taking into account the number of patients who are at risk of the transition and the patients who finally have done the transition. For each row of each matrix the number of patients at risk will be different.

We present here the estimation of the matrices $\mathbf{{P}_{(1)}}$ and $\mathbf{{P}_{(2)}}$. The estimation of the rest of the matrices  is similar, except for matrices $\mathbf{{P}_{(6)}}$ and $\mathbf{{P}_{(7)}}$ ($6$ and $7$ are absorbent states) which are null matrices except for the elements $(6,6)$ and $(7,7)$ which are equal to 1.

In order to estimate the matrix $\mathbf{{P}_{(1)}}$ we start from all patients who have been hospitalized with entry in state NSP. The day after a patient has been hospitalized he/she  can still be at NSP or can move to SP, Discharge or Death. So rows 1, 2, 6 and 7 are the only ones with probabilities different from 0. In order to estimate the probability cells in row 1 of $\mathbf{{P}_{(1)}}$ we consider, for each time $s$, all the patients who have been at least two consecutive days in NSP, that is,  $\sum_{s=2}^{43} Y_{11}(s-1)=10577$ (note here that 43 is the maximum number of days a patient has been two consecutive days in NSP). From those $10577$ patients at risk,  the number of patients who have stayed in NSP the  next day is $\sum_{s=2}^{43}  N_{111}(s)=8919$, while $\sum_{s=2}^{43}  N_{112}(s)=257$ have transited to  SP, $\sum_{s=2}^{43}  N_{116}(s)=1369$ have been discharged and, finally, $\sum_{s=2}^{43}  N_{117}(s)=32$ have died.  We proceed analogously  for the estimation  of  the probability cells in row 2 of $\mathbf{{P}_{(1)}}$ starting  with those  $\sum_{s=2}^{36}  Y_{12}(s-1)=411$ patients who have moved  to SP from NSP the next day.  
  
\[\mathbf{\widehat{P}_{(1)}}=\begin{pmatrix}
        \frac{8919}{10577} &  \frac{257}{10577} & 0 & 0 & 0 &  \frac{1369}{10577} &  \frac{32}{10577}\\[4pt]
        0 & \frac{253}{411} & \frac{3}{411} & \frac{92}{411} & \frac{62}{411} & 0 & \frac{1}{411}\\[4pt]
        0 & 0 & 0 & 0 & 0 & 0 & 0\\
        0 & 0 & 0 & 0 & 0 & 0 & 0\\
        0 & 0 & 0 & 0 & 0 & 0 & 0\\
        0 & 0 & 0 & 0 & 0 & 1 & 0\\
        0 & 0 & 0 & 0 & 0 & 0 & 1\\
      \end{pmatrix}\]

 \vskip 15pt     
In a similar way we estimate $\mathbf{{P}_{(2)}}$. In this case we start with all patients who have been at state SP at any time. The next day these patients can still be at SP or can move to Recovery, NIMV, IMV or Death. So rows 1 and 6 are 0 because there is no direct transition from SP to NSP nor to Discharge. For row 2, the number of patients at risk, that is, the number of patients spending two consecutive times in state SP is $\sum_{s=2}^{47}  Y_{22}(s-1)=2668$. Row 3  starts with those patients who have moved from  SP to Recovery, a total of $\sum_{s=2}^{47} Y_{23}(s-1)=223$. Analogously for rows 4 and 5, $\sum_{s=2}^{14} Y_{24}(s-1)=214$  patients have transited immediately from  SP to NIMV while $\sum_{s=2}^{23} Y_{25}(s-1)=166$ patients moved from SP to  IMV. Probability matrices $\mathbf{{P}_{(3)}}$, $\mathbf{{P}_{(4)}}$ and $\mathbf{{P}_{(5)}}$ are estimated proceeding in an analogous manner, each one starting from patients in states Recov, NIMV and IMV, respectively.

\vskip 15pt
  \[\mathbf{\widehat{P}_{(2)}}= \begin{pmatrix}
   0 & 0 & 0 & 0 & 0 & 0 & 0\\
  0 & {\frac{2307}{2668}} &{\frac{220}{2668}} & {\frac{68}{2668}} & {\frac{49}{2668}} & 0 &{\frac{24}{2668}}\\[4pt]
  0 & 0 & {\frac{207}{223}}& 0 & 0 & {\frac{16}{223}} & {0}\\[4pt]
0 & 0 & {\frac{6}{214}} & {\frac{159}{214}}&{\frac{45}{214}} & 0 &{\frac{4}{214}}\\[4pt]
    0 & 0 & {{\frac{3}{166}}} & 0 & {{\frac{160}{166}}} & 0 & {{\frac{3}{166}}}\\
      0 & 0 & 0 & 0 & 0 & 0 & 0\\
        0 & 0 & 0 & 0 & 0 & 0 & 1
    \end{pmatrix}\]

%%%%%%%%%%%%%%%%%%%%%%%%%%%%%%%%%%%%%%%%%%%%%%%
\subsection{Prediction via Chapman-Kolmogorov equations}
%%%%%%%%%%%%%%%%%%%%%%%%%%%%%%%%%%%%%%%%%%%%%%%
The Markov test computed in Section \ref{sec:MarkovTest} rejects the first order Markov assumption for three of the  four
transitions from Severe Pneumonia (SP):  to Non Invasive Mechanical Ventilation (NIMV), Invasive Mechanical Ventilation (IMV) and Recovery (Recov), indicating that whether or not the patient was diagnosed with Non Severe Pneumonia (NSP) marks a difference in his/her  prognosis.
A second order model allows the prediction of  the time to future events as a function of the diagnostic when they were hospitalized. Chapman-Kolomogorov extension in Theorem \ref{eq:C-K_M} is the key to the corresponding probabilities.

\begin{figure}[H] 
        \centering
        \includegraphics[scale=0.4]{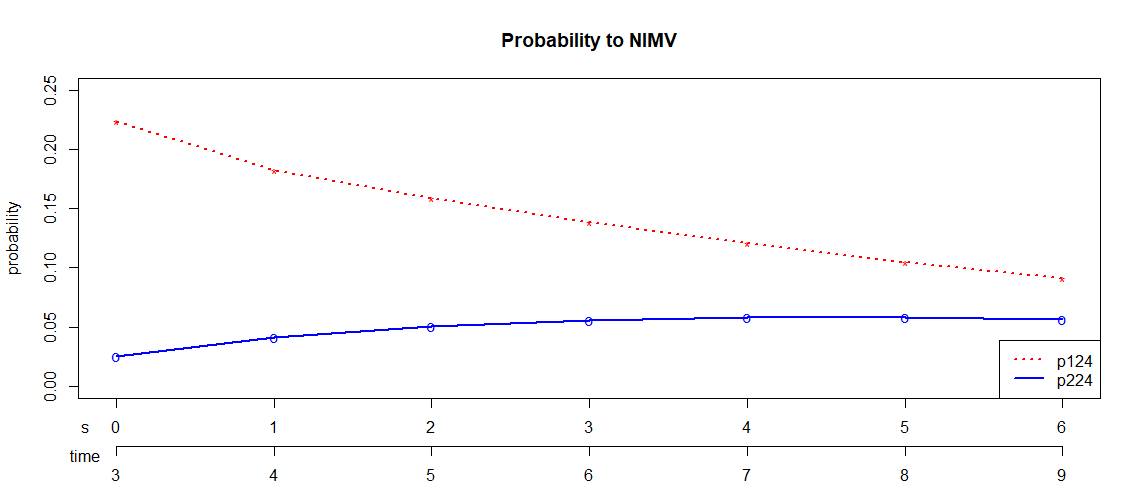}
        \caption{Probabilities from SP to NIMV for patients who had been 2 consecutive days in SP (line) compared with patients who have been one day in NSP and one day in SP (dots). $p_{224}=P(X_{3+s}=4|X_2=2,X_1=2)$ vs $p_{124}=P(X_{3+s}=4|X_2=2,X_1=1)$}
        \label{img:224-124}
    \end{figure}

For the transition from SP $\rightarrow$ NIMV, we will compute for $s\in\{0,\ldots,6\}$ the probabilities $p_{224}=P(X_{3+s}=4|X_2=2,X_1=2)$, that is,  the probability to NIMV for patients who arrive at the hospital with a SP diagnosis and they still were in SP the second day. And also $p_{124}=P(X_{3+s}=4|X_2=2,X_1=1)$ the probability to NIMV for patients with NSP at admission  who had moved to SP the second day. We have  plotted these probabilities in Figure \ref{img:224-124}, where we can we see  how important is the initial state for the initial times. The probability of moving to NIMV of patients initially diagnosed with SP ($X_1=2, X_2=2$) is  much smaller than the probability of moving to NIMV of patients initially diagnosed with NSP ($X_1=1, X_2=2$). These two probabilities close the gap as days go by.

\begin{figure}[H]
        \centering
        \includegraphics[scale=0.4]{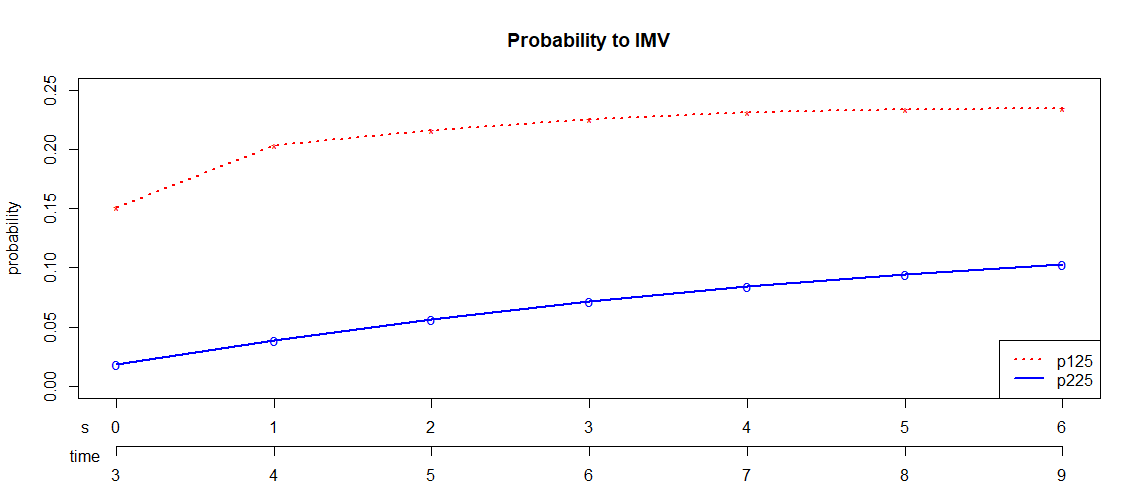}
        \caption{Probabilities from SP to IMV for patients who had been 2 consecutive days in SP (line) compared with patients who have been one day in NSP and one day in SP (dots). $p_{225}=P(X_{3+s}=5|X_2=2,X_1=2)$ vs $p_{125}=P(X_{3+s}=5|X_2=2,X_1=1)$} 
        \label{img:225-125}
    \end{figure}

The same  type of plot is depicted in Figure \ref{img:225-125} to  study the transition SP $\rightarrow$ IMV. In this case the patients  are also splitted based on their initial state: NSP or SP. As in Figure \ref{img:224-124}, the probability of moving to IMV of patients initially diagnosed with SP ($X_1=2, X_2=2$) is  much smaller than the probability of moving to IMV of patients initially diagnosed with NSP ($X_1=1, X_2=2$). However, both  probabilities increase over the time and their difference is kept along the next days. This reveals the different prognosis for needing respiratory mechanical ventilation (IMV)  among those  patients initially diagnosed with NSP versus being diagnosed with SP. 

%%%%%%%%%%%%%%%%%%%%%%%%%%%%%%%%%%%%%%%%%%%%%%%%%%%%%%%%%%%
%%%%%%%%%%%%%%%%%%%%%%%%%%%%%%%%%%%%%%%%%%%%%%%%%%%%%%%%%%%
\section{Discussion}\label{discussion}
%%%%%%%%%%%%%%%%%%%%%%%%%%%%%%%%%%%%%%%%%%%%%%%%%%%%%%%%%%%
%%%%%%%%%%%%%%%%%%%%%%%%%%%%%%%%%%%%%%%%%%%%%%%%%%%%%%%%%%%

In this paper we have introduced a second order Markov multistate model and we have developed an extension of  the Chapman-Kolmogorov equations to compute $r$-step transition probabilities. We have used the DIVINE COVID-19 data to estimate the transition probabilities and to predict probabilities to NIMV and  IMV in terms of the states where a patient was during the first 2 days of his/her hospitalization.

It should be mentioned that Chapman-Kolmogorov extension is  based on a discretization of the time scale and is only computed conditionally to two-consecutive times.  The computation of the transition probabilities  if the previous two times were not consecutive would be unfeasible because
the window  between two not consecutive times might include a subject in several states.

 In this paper we sketch two different ways to estimate  the transition probabilities. The first one  using the Bernoulli probabilities, which is the one used to compute the transition probabilities in the COVID illustration example and the second using  conditional probabilities. Since the data from the DIVINE project was collected one year  after the end of the first wave, we have complete registries  and, for this reason, we have so far only developed both methods for complete (uncensored) data. Nevertheless, it is indeed relevant to extend these estimators to account  
 for right-censored data.
The second method of estimation  presented in Subsection \ref{sec:EstCondProb}  gives a clue of how we could proceed to account  
 for right-censored data. This estimator, an average of the  ratios, for each time, of those subjects doing an specific transition among the number of subjects at risk, has an analogy to the Nelson-Aalen estimator for the cumulative hazard function. For a thorough statistical analysis, the derivation of the variance of these estimators as well as of their 
 asymptotic distribution is needed. Furthermore, estimators for the state occupation probabilities and  for the transition intensities for complete and right-censored data are as well a topic of interest. All these ideas remain open for our future research.

A second order Markov model could had been transformed into a first order Markov model by redefining the state space. This would be possible creating extra states formed by direct 1-step transitions. For instance, instead of one unique Death state we could have defined 3 new states formed by those patients arriving to Death from NIMV, IMV or Recovery. The advantage of these  new states is clear because  we would be able
to apply all the knowledge on first order Markov models. However, the number of states and transitions of the new model will increase substantially and the interpretation will become cumbersome. Furthermore, since the number of parameters to estimate will increase and,  so does the needed sample sample size 
to estimate all of them, the second order Markov approach is preferable.
%A second order Markov model solves this issue without losing the idea of dependence on the two previous states.

\section*{Acknowledgements} This research has been funded by the Ministerio de Ciencia e Innovaci\'on (Spain) [PID2019-104830RB-I00/ DOI(AEI): 10.13039/501100011033] and by Generalitat de Catalunya  through the projects 2020PANDE00148 and 01421 SGR-Cat 2021. We are indebted to our colleagues in the DIVINE group  for their clever contributions and dedicated time.

\bibliographystyle{chicago}
\bibliography{sn-bibliography}

\begin{thebibliography}{}

\bibitem[\protect\citeauthoryear{Chakladar, Liao, Landau, Gamalo, and
  Wang}{Chakladar et~al.}{2022}]{Chakladar}
Chakladar, S., R.~Liao, W.~Landau, M.~Gamalo, and Y.~Wang (2022).
\newblock {Discrete Time Multistate Model With Regime Switching for Modeling
  COVID-19 Disease Progression and Clinical Outcomes}.
\newblock {\em Statistics in Biopharmaceutical Research\/}~{\em 14\/}(1),
  52--66.

\bibitem[\protect\citeauthoryear{Ching, Fung, and Ng}{Ching
  et~al.}{2003}]{Ching(2003)}
Ching, W.~K., E.~S. Fung, and M.~K. Ng (2003).
\newblock {A higher-order Markov model for the Newsboy's problem}.
\newblock {\em Journal of the Operational Research Society\/}~{\em 54\/}(3),
  291--298.

\bibitem[\protect\citeauthoryear{de~Wreede, Fiocco, and Putter}{de~Wreede
  et~al.}{2011}]{mstate}
de~Wreede, L., M.~Fiocco, and H.~Putter (2011).
\newblock {\texttt{mstate}: An R Package for the Analysis of Competing Risks
  and Multi-State Models}.
\newblock {\em Journal of Statistical Software\/}~{\em 38\/}(7), 1--30.

\bibitem[\protect\citeauthoryear{Garmendia, Cort\'es, and
  G\'omez~Melis}{Garmendia et~al.}{2023}]{Garmendia}
Garmendia, L., J.~Cort\'es, and G.~G\'omez~Melis (2023).
\newblock {MSMpred: Interactive modelling and prediction of individual
  evolution via multistate models}.
\newblock {\em {(submitted)}\/}.

\bibitem[\protect\citeauthoryear{Hougaard}{Hougaard}{1999}]{Hougaard}
Hougaard, P. (1999).
\newblock {Multi-state models: a review}.
\newblock {\em Lifetime Data Anal.\/}~{\em 5\/}(3), 239--264.

\bibitem[\protect\citeauthoryear{Islam and Chowdhury}{Islam and
  Chowdhury}{2006}]{Islam(2006)}
Islam, M.~A. and R.~I. Chowdhury (2006).
\newblock {A higher order Markov model for analyzing covariate dependence}.
\newblock {\em Applied Mathematical Modelling\/}~{\em 30\/}(6), 477--488.

\bibitem[\protect\citeauthoryear{Kay}{Kay}{1986}]{Kay}
Kay, R. (1986).
\newblock {A Markov Model for Analysing Cancer Markers and Disease States in
  Survival Studies}.
\newblock {\em Biometrics\/}~{\em 42\/}(4), 855--865.

\bibitem[\protect\citeauthoryear{Lin, Wei, and Ying}{Lin
  et~al.}{1993}]{Lin(1993)}
Lin, D.~Y., L.~J. Wei, and Z.~Ying (1993).
\newblock Checking the cox model with cumulative sums of martingale-based
  residuals.
\newblock {\em Biometrika\/}~{\em 80\/}(3), 557--572.

\bibitem[\protect\citeauthoryear{Logan}{Logan}{1981}]{Logan(1981)}
Logan, J.~A. (1981).
\newblock {A structural model of the higher-order Markov process incorporating
  reversion effects}.
\newblock {\em The Journal of Mathematical Sociology\/}~{\em 8\/}(1), 75--89.

\bibitem[\protect\citeauthoryear{Pallar\`es, Teb\'e, Abelenda-Alonso, Rombauts,
  Oriol, Simonetti, Rodr\'iguez-Molinero, Izquierdo, D\'iaz-Brito, Molist,
  G\'omez~Melis, Carratal\`a, Videla, and study groups}{Pallar\`es
  et~al.}{2023}]{Pallares}
Pallar\`es, N., C.~Teb\'e, G.~Abelenda-Alonso, A.~Rombauts, I.~Oriol, A.~F.
  Simonetti, A.~Rodr\'iguez-Molinero, E.~Izquierdo, V.~D\'iaz-Brito, G.~Molist,
  G.~G\'omez~Melis, J.~Carratal\`a, S.~Videla, and M.~study groups (2023).
\newblock {Characteristics and Outcomes by Ceiling of Care of Subjects
  Hospitalized with COVID-19 During Four Waves of the Pandemic in a
  Metropolitan Area: A Multicenter Cohort Study}.
\newblock {\em Infectious diseases and therapy\/}~{\em 12\/}(1), 273--289.

\bibitem[\protect\citeauthoryear{Piulachs, Langhor, Besal\'u, Pallar\`es,
  Carratal\`a, Teb\'e, and G\'omez~Melis}{Piulachs et~al.}{2023}]{Piulachs}
Piulachs, X., K.~Langhor, M.~Besal\'u, N.~Pallar\`es, J.~Carratal\`a,
  C.~Teb\'e, and G.~G\'omez~Melis (2023).
\newblock {Semi-Markov multistate approaches for multicohort event history
  data}.
\newblock {\em (submitted)\/}.

\bibitem[\protect\citeauthoryear{Raftery}{Raftery}{1985}]{Raftery(1985)}
Raftery, A.~E. (1985).
\newblock {A Model for High-Order Markov Chains}.
\newblock {\em Journal of the Royal Statistical Society: Series B
  (Methodological)\/}~{\em 47\/}(3), 528--539.

\bibitem[\protect\citeauthoryear{Rodr\'iguez-Girondo and de~U\~{n}a
  \'Alvarez}{Rodr\'iguez-Girondo and de~U\~{n}a
  \'Alvarez}{2012}]{Rodriguez-Girondo}
Rodr\'iguez-Girondo, M. and J.~de~U\~{n}a \'Alvarez (2012).
\newblock {A nonparametric test for Markovianity in the illness-death model}.
\newblock {\em Statistics in Medicine\/}~{\em 31\/}(30), 4416--4427.

\bibitem[\protect\citeauthoryear{Shamshad, Bawadi, {Wan Hussin}, Majid, and
  Sanusi}{Shamshad et~al.}{2005}]{Shamshad}
Shamshad, A., M.~Bawadi, W.~{Wan Hussin}, T.~Majid, and S.~Sanusi (2005).
\newblock First and second order markov chain models for synthetic generation
  of wind speed time series.
\newblock {\em Energy\/}~{\em 30\/}(5), 693--708.

\bibitem[\protect\citeauthoryear{Shorrocks}{Shorrocks}{1976}]{Shorrocks(1976)}
Shorrocks, A.~F. (1976).
\newblock {Income Mobility and the Markov Assumption}.
\newblock {\em The Economic Journal\/}~{\em 86\/}(343), 566--578.

\bibitem[\protect\citeauthoryear{Titman and Putter}{Titman and
  Putter}{2020}]{TitmanPutter}
Titman, A.~C. and H.~Putter (2020).
\newblock {General tests of the Markov property in multi-state models}.
\newblock {\em Biostatistics\/}~{\em 23\/}(2), 380--396.

\bibitem[\protect\citeauthoryear{Tong}{Tong}{1975}]{Tong(1975)}
Tong, H. (1975).
\newblock {Determination of the order of a Markov chain by Akaike's information
  criterion}.
\newblock {\em Journal of Applied Probability\/}~{\em 12\/}(3), 488--497.

\end{thebibliography}
\end{document}